\documentclass[12pt]{amsart}

\newtheorem{theorem}{Theorem}[section]

\theoremstyle{definition}
\newtheorem{definition}[theorem]{Definition}
\newtheorem{example}[theorem]{Example}

\newtheorem{corollary}[theorem]{Corollary}
\newtheorem{proposition}[theorem]{Proposition}

\theoremstyle{remark}
\newtheorem{remark}[theorem]{Remark}

\numberwithin{equation}{section}
\numberwithin{theorem}{section}

\newcommand\bel[1]{\begin{equation}\label{#1}}
\newcommand\ee{\end{equation}}

\newcommand\und[1]{{\tt{#1}}}

\newcommand{\ol}[1]{{\overline{#1}}}
\newcommand{\wt}[1]{{\widetilde{#1}}}

\newcommand{\bR}{{\mathbb{R}}}
\newcommand{\bE}{{\mathbb{E}}}

\newcommand{\bA}{{\mathbf{V}}}

\newcommand{\bH}{{\mathbf{H}}}
\newcommand{\bfX}{{\mathbf{X}}}

\newcommand{\cA}{{\mathcal{A}}}
\newcommand{\cB}{{\mathcal{B}}}
\newcommand{\cF}{{\mathcal{F}}}
\newcommand{\cI}{{\mathcal{I}}}
\newcommand{\cJ}{{\mathcal{J}}}
\newcommand{\cK}{{\mathcal{K}}}
\newcommand{\cP}{{\mathcal{P}}}
\newcommand{\cR}{{\mathcal{R}}}

\newcommand{\mfX}{{\mathfrak{X}}}
\newcommand{\mfB}{{\mathfrak{B}}}

\begin{document}

\title[Chaos solutions]
{Solving SPDE{s} driven by colored noise:
a chaos approach}

\author{S. V. Lototsky\ }
 \address[S. V. Lototsky]{Department of
Mathematics, USC\\ Los Angeles, CA 90089}
 \email[S. V. Lototsky]{lototsky@math.usc.edu}
\urladdr{http://www-rcf.usc.edu/$\sim$lototsky}
\thanks{S. V. Lototsky  acknowledges  support
 from the NSF CAREER award DMS-0237724}

 \author{\ K. Stemmann}
\address[K. Stemmann]{Department of
Mathematics, USC\\ Los Angeles,
CA 90089}
 \thanks{
The work of K. Stemmann was partially supported by
the NSF Grant DMS-0237724}

\begin{abstract}
An It\^{o}-Skorokhod  bi-linear equation
driven by infinitely many
independent colored noises is considered in
a normal triple of Hilbert spaces. The special
feature of the equation is the
appearance of the Wick product in the definition
of the It\^{o}-Skorokhod integral, requiring
innovative approaches to computing the
solution.    A chaos
expansion of the solution is derived and
several truncations of this expansion are
studied. A recursive approximation of the
solution is suggested and the corresponding
approximation error bound is computed.
\end{abstract}

\maketitle

\section{Introduction}
\setcounter{equation}{0}

Stochastic differential equations driven
by Gaussian white noise
 are well-studied; see, for example,
the book \cite{Oksendal} for ordinary
differential equation
and the book \cite{Roz} for equations with
partial derivatives. The underlying stochastic process
in these equations is the standard Brownian motion
$W$, which is a square-integrable Gaussian martingale
with continuous trajectories and independent
increments. A lot less is known about
equations driven by {\em colored noise}, when the
underlying process is still Gaussian, but
no longer has independent increments.
 An important example
is the fractional Brownian motion $W^{H}$,
$H\in (0,1)$,  which
coincides with the standard Brownian motion $W$ for
 $H=1/2$ and is not a
semi-martingale for all $H\not=1/2$.
It is the lack of the semi-martingale property
that makes the analysis difficult at the
very basic level, the definition of the
corresponding stochastic integral.
Several versions of the
stochastic integral with respect to $W^H$ have been
proposed \cite{AMN1, DH, DU, DHP, KKA, Lin, Za}.
Unlike the standard Brownian motion,
 different approaches  such as It\^{o}-type vs.
 Stratonovich-type
integral or  path-vise vs. mean-square definition,
 become much more
difficult to reconcile. The paper by V. Pipiras
 and M. Taqqu
\cite{PT} describes  the main difference between
integration with
respect to $W$ and $W^H$ for $H\not=1/2$, and
 demonstrates the
resulting technical difficulty. Beside a purely
theoretical
interest, fractional Brownian motion appears to
be a natural
replacement of the standard Brownian motion in
certain applied problems
\cite{A, AC, BP, LeB, HO, KBR}. In many
such problems, it is possible to avoid most of
the technical issues
related to the definition of the stochastic integral
  by considering additive noise.

One way to streamline the analysis of differential
equations driven by multiplicative
 noise (also known as {\tt bilinear equations})
  is to use the Wick product.
It is known \cite{HOUZ} that, with standard
Brownian motion,
Wick product and usual calculus lead to the same
results as the
usual product and the It\^{o} calculus.
While the use of Wick
product has been questioned as a
modeling tool for certain
applications in economics and
finance \cite{BH}, it is still an effective
tool for theoretical investigations, corresponding
to the It\^{o}-Skorokhod integral in the white
noise analysis.

Successful study of  differential
equations driven by multiplicative
colored noise also requires a convenient representation
of the underling Gaussian process.
Traditionally, a Gaussian
process is defined by its mean and covariance
functions, but then the
definition of the integral immediately leads
to a number of technical
 conditions on these functions \cite{AMN}.
 An alternative definition is possible \cite{LtSt}
 by combining the ideas from
the theory of generalized Gaussian fields
  \cite{GV, Maj},
  the white noise theory \cite{HKPS, HOUZ},
  and the Malliavin calculus
  \cite{Mal, Nualart}. This approach to stochastic
  integration is  used in this paper
and is outlined in Section \ref{Sec:CN}

Numerical methods for stochastic ordinary
differential equations driven
by white noise are a well-developed
subject \cite{KlPl, Mil}. Equations with
partial derivatives have been mostly studied in
connection with optimal nonlinear
 filtering
 \cite{Bennaton, BGR, GermPicc, KIto, ItoRoz}.
 In the end, all
these numerical methods have obvious counterparts
in the numerical analysis
of the deterministic equations (Galerkin method,
Euler method, finite difference
and finite element methods, operator splitting
method, etc.)

The main difficulty in the numerical analysis of
equations driven by colored noise is the use of
the Wick product: unlike the usual
product, Wick product is not an operation readily
performed by a computer.
Since the Wick product is relatively easy to compute
for Hermite
polynomials of Gaussian random variables, an
 implementation of this operation should be based
 on the chaos expansion, and then
 a truncation of this expansion becomes
a natural numerical approximation of the solution.

While  a truncated
chaos expansion has been investigated before
as an approximation of the solution of a stochastic
equation,
in particular, for equations of optimal nonlinear
filtering \cite[etc.]{BdKl, LMR},
 this approximation was always just another
 possibility of solving the
equation numerically. By contrast, for equations
 driven by colored noise,
chaos expansion appears to be the only possibility to
compute the solution. The
 chaos expansion of the solution of a bi-linear
equation driven by infinitely many independent
colored noises is studied below
in Section \ref{Sec:solS};
various truncations of this expansion are studied
 in Section \ref{Sec:appr}.

\section{Colored Noise}
\label{Sec:CN}

Let $(\Omega, \cF, \mathbb{P})$ be a complete
probability space.

\begin{definition}
\label{def:field}
(a) A \und{colored noise} $\mfX$ on $L_2((0,T))$
 with the covariance operator $\cR$ is a collection
 of zero-mean Gaussian random
variables $\mfX(f),\ f\in L_2((0,T)),$ with the
property
\bel{eq:1.1}
\bE\Big( \mfX(f) \mfX(g) \Big)=
\int_0^T (\cR f)(s)g(s)ds,\ f,\,g\in L_2((0,T)),
\ee
where $\cR$ is a
bounded linear operator on $L_2((0,T))$.

(b) A \und{representation operator} of $\mfX$ is a
bounded linear operator $\cK$ on $L_2((0,T))$
such that  $\cK\cK^*=\cR$, where $\cK^*$ is the
adjoint of $\cK$.
\end{definition}

\begin{remark}
\label{rem-rep}
(a) The white noise on $L_2((0,T))$
corresponds to $\cR=I$,
the identity operator \cite{HKPS}.

(b) Since $\cR$ is a self-adjoint non-negative
operator, a representation operator always exists.

\end{remark}

For $t>0$ let
$$
\chi_t(s)=
\begin{cases}
1,& 0\leq  s\leq t,\\
0,& {\rm otherwise,}
\end{cases}
$$
denote the characteristic function of the
interval $[0,t]$.

\begin{example}
If $\mfX$ is white noise, then
 direct computations show that
$W(t)=\mfX(\chi_t)$ is a standard Brownian motion and
\bel{eq:wn-rep}
\mfX(f)=\int_0^T f(s)dW(s),\ f\in L_2((0,T)).
\ee
\end{example}

\begin{example}
\label{exe:OU}
{\bf (Ornstein-Uhlenbeck Noise)}

For $b>0$ and $t,s\in [0,T]$ define
$$
K(t,s)=-be^{-b(t-s)}.
$$
Then the operator
$$
\cK: f(t) \mapsto (\cK f)(t)=f(t)+\int_0^t K(t,s)f(s)ds
$$
is bounded on $L_2((0,T))$:
$$
\int_0^T |(\cK f)(s)|^2ds \leq (1+\sqrt{bT})^2
\int_0^T |f(s)|^2ds.
$$
If $\mfX$ is the colored noise corresponding
to the representation operator $\cK$, then
a straightforward computation shows that
$\mfX(\chi_t),\ t\in (0,T),$ is the
Ornstein-Uhlenbeck process with covariance
$$
\bE\big(\mfX(\chi_t)\mfX(\chi_s)\big)=
e^{-b|t-s|}.
$$
 Accordingly, we call $\mfX$
the \und{Ornstein-Uhlenbeck noise} with parameter $b$.
\end{example}

\begin{example}
\label{Exe-fBm}
{\bf (Fractional White Noise)}

 For $H\in (1/2,1)$ and $t,s \in [0,T]$ define
$$
K(t,s)=C_H\left(H-\frac{1}{2}\right)
 \left(\frac{t}{s}\right)^{\frac{1}{2}-H}
 (t-s)^{H-\frac{3}{2}}\,\chi_t(s),
$$
where
$$
C_H=\left(\frac{2H\Gamma\left(\frac{3}{2}-H\right)}
{\Gamma\left(H+\frac{1}{2}\right)
\Gamma(2-2H)}\right)^{\frac{1}{2}}
$$
and $\Gamma$ is the Gamma function.
Then the operator
$$
 \cK: f(t)\mapsto (\cK f)(t)=\int_0^t K(t,s) f(s)ds
$$
is bounded on $L_2((0,T))$ \cite{LtSt,Nualart2}:
\begin{equation}
\label{eq:fbm-norm}
\int_0^T |(\cK f)(s)|^2ds \leq
\frac{H(2H-1)\,\Gamma\left(H-\frac{1}{2}\right)}
{\Gamma\left(H+\frac{1}{2}\right)}\ T^{2H-1}
\int_0^T |f(s)|^2ds.
\end{equation}
If $\mfX$ is the colored noise corresponding
to the representation operator $\cK$, then
$\mfX(\chi_t),\ t\in (0,T),$ is the
fractional Brownian motion with Hurst parameter $H$
\cite{Nualart2}.
Accordingly, we call $\mfX$
the \und{$H$-fractional white noise}.
\end{example}


In general, if $\mfX$ is a colored noise on
$L_2((0,T))$, then $X(t)=\mfX(\chi_t)$ is a zero-mean
Gaussian process. Thus,  one can interpret
$\mfX$ as a collection of integrals $\int_0^Tf(s)dX(s)$
for {\em deterministic} $f\in L_2((0,T))$. Our
study of bi-linear equations (linear equations
with multiplicative noise) requires an extension
of $\mfX$ to {\em random} $f$ and is based on
the following generalization of \eqref{eq:wn-rep}.

\begin{proposition}
\label{prop:rep}
For every colored noise $\mfX$ on $L_2((0,T))$
with a representation operator $\cK$, there
exists a unique standard Brownian motion $W=W(t)$
such that
\bel{eq:rep000}
\mfX(f)=\int_0^T (\cK^* f)(s)dW(s),\ f\in L_2((0,T)).
\ee
\end{proposition}

\begin{proof}
Relation \eqref{eq:rep000} certainly defines a colored
noise on $L_2((0,T))$; we omit the technical proof that
the corresponding Brownian motion can be found for
every colored noise $\mfX$ \cite{LtSt}.
\end{proof}

\begin{definition}
\label{def-rep1}
A pair $(\cK,W)$, where $\cK$ is a bounded linear
operator on $L_2((0,T))$ and $W$ is a standard
Brownian motion, is called a \und{representation}
of the colored noise $\mfX$ if \eqref{eq:rep000}
holds.
\end{definition}

For random $f$, we now define $\mfX(f)$ according to
\eqref{eq:rep000}, where the stochastic integral
is understood in the \und{It\^{o}-Skorokhod sense}
\cite{Nualart}. A equivalent, but less convenient,
definition of $\mfX(f)$ for random $f$ is
possible in intrinsic terms without using a
representation of $\mfX$ \cite{LtSt}.

\section{Chaos Solution: Existence and Regularity}
\label{Sec:solS}

Let $(\Omega, \cF, \mathbb{P})$ be a complete
probability space and $\mfX_{\ell},\ \ell\geq 1$,
a collection of independent colored noises on
$L_2((0,T))$.

In this section we study the equation
 \bel{eq:LSp-eqm1}
u(t)=u_0+\int_0^t (\cA u(s)+F(s)) ds
+ \sum_{\ell=1}^{\infty}
\mfX_{\ell}(\chi_t\cB_{\ell} u+\chi_tG_{\ell}),
\ee
$t\in [0,T]$, in a normal triple $(\bfX,\bH,\bfX')$
 of Hilbert spaces. In particular, we assume that
 $\cA$ and each $\cB_{\ell}$ are
 bounded linear operators from $\bfX$ to
 $\bfX'$, $u_0\in L_2(\Omega; \bfX')$,
 $F, G_{\ell} \in L_2(\Omega;L_2((0,T);\bfX'))$.

 \begin{remark} While it is tempting to re-write
\eqref{eq:LSp-eqm1} as
\bel{eq:LSp-eqm-intuit}
u(t)=u_0+\int_0^t (\cA u(s)+F(s)) ds
+ \sum_{\ell=1}^{\infty}\int_0^t
(\cB_{\ell} u(s)+G_{\ell}(s))dX_{\ell}(s),
\ee
where $X_{\ell}(t)=\mfX_{\ell}(\chi_t)$,
we will use a more rigorous  form \eqref{eq:LSp-eqm1}.
\end{remark}

  By analogy with equations
 driven by white noise \cite{Roz}, we define a
 \und{variational solution} of \eqref{eq:LSp-eqm1}
 as a random element $u$ with values in
 $L_2(\Omega; L_2((0,T);\bfX))$ such that,
 for every $v\in \bfX$, the equality
 \begin{equation}
 ( u(t), v)_{\bH}
 =\langle u_0, v\rangle
 +\int_0^t \langle (\cA u(s)+F(s)), v\rangle ds
+ \sum_{\ell=1}^{\infty}
\mfX_{\ell}(\chi_t
\langle \cB_{\ell} u+G_{\ell}, v\rangle)
\end{equation}
holds in $\bfX'$ on the same set of probability one
for all $t\in [0,T]$; by $\langle \cdot, \cdot \rangle$
we denote the duality between $\bfX'$ and $\bfX$
relative to the inner product
$(\cdot, \cdot)_{\bH}$  of $\bH$.

Unfortunately, the
current development of the colored noise calculus
is not sufficient to establish existence of
a variational solution of \eqref{eq:LSp-eqm1}.
Accordingly, we introduce a weaker notion of
solution, called a {\em chaos solution},
using a Fourier series expansion in the space of
square integrable random processes.

We start with some auxiliary constructions.
Let $(\cK_{\ell}, W_{\ell})$ be  a representation
of the colored noise $\mfX_{\ell}$ and let
$\{h_k,\, k\geq 1\}$ be an orthonormal basis
in $L_2((0,T))$. Define random variables
\bel{eq:xi-k}
\xi_{k\ell}=\int_0^T h_k(t)dW_{\ell}(t).
\ee
Let $\cJ$ be the collection of multi-indices
 $\alpha=\{\alpha_{k\ell},\ k,\ell\geq 1\}$.
Each $\alpha\in \cJ$ has  non-negative integer elements
$\alpha_{k\ell}$ and
$\sum_{k,\ell}\alpha_{k\ell}<\infty$.
We will use the notations
\begin{equation}
|\alpha|=\sum_{k,\ell}\alpha_{k\ell},\
\alpha!=\prod_{k,\ell}\alpha_{k\ell}!.
\end{equation}
By $(0)$ we denote the multi-index $\alpha$
 with $|\alpha|=0$ and by $\epsilon_{ij}$,
 the multi-index $\alpha$ with $|\alpha|=1$ and
 $\alpha_{ij}=1$.

 \begin{remark} If there is only one colored noise,
 then the entries of $\alpha$ have only one
 index: $\alpha=\{\alpha_k,\ k\geq 1\}$.
 \end{remark}

For $\alpha\in \cJ$ define
\bel{eq:xi-al}
\xi_{\alpha}=
\prod_{k,\ell\geq 1}
\frac{H_{\alpha_{k\ell}}(\xi_{k\ell})}
{\sqrt{\alpha_{k\ell}!}},
 \ee
where, for an integer $n\geq 0$,
   $H_n=H_n(t)$ is the $n$-th \und{Hermite polynomial}
\begin{equation}
\label{eq:HerPol}
 H_{n}(t)=(-1)^ne^{t^{2}/2}
\frac{d^{n}}{dt^{n}}e^{-t^{2}/2}.
\end{equation}

Recall the definition of the Wick product \cite{HOUZ}:
\bel{eq:wick}
 H_m(\xi_{ij})\diamond H_n(\xi_{k\ell})=
 \begin{cases}
 H_{m+n}(\xi_{ij}),& {\rm if \ } i=k
          {\rm \ and \ }j=\ell\\
H_m(\xi_{ij})H_{n}(\xi_{k\ell}),& {\rm otherwise}.
\end{cases}
\ee
In particular, for $m\geq 1$,
\begin{equation}
H_{m}(\xi_{k\ell})
=\underset{m\,\mathrm{times}}
{\underbrace{\xi_{k\ell}\diamond
\cdots\diamond\xi_{k\ell}}}.
\end{equation}
Then every multi-index $\alpha$ is uniquely
characterized by the collection \\
$\{(k_1, \ell_1), \ldots, (k_n,\ell_n)\}$,
called the \und{characteristic set} of $\alpha$,
such  that $k_1\leq k_2\leq\cdots\leq k_n$,
$\ell_i\leq \ell_{i+1}$ if $k_i=k_{i+1}$, and
\begin{equation}
\label{xial-alt}
\xi_{\alpha}= \frac{\xi_{k_{1}\ell_1}
\diamond\xi_{k_{2}\ell_2}
\diamond\cdots\diamond\xi_{k_{n}\ell_n}}
{\sqrt{\alpha!}}.
\end{equation}

\begin{proposition}
\label{prob:stoch-bas}

(a) The collection $\{\xi_{\alpha},\ \alpha\in \cJ\}$
is an orthonormal basis in the space of square
integrable random variables that are measurable
with respect to the $\sigma$-algebra $\cF^W$ generated by
the Brownian motions $W_{\ell}$, $\ell\geq 1$,
on $[0,T]$.

(b) Let  $\eta$ be a square-integrable
 $\cF^W$-measurable random
element with values in $L_2((0,T))$. Define
\begin{equation}
{\eta}_{\alpha}(t)=\bE ({\eta}_k(t)\xi_{\alpha}).
\end{equation}
 Then
\bel{eq:int-chaos}
\mfX_{\ell}(\eta)=\sum_{\alpha\in \cJ}
\left(\sum_{k\geq 1}\sqrt{\alpha_{k\ell}}\,
\int_0^T{\eta}_{\alpha-\epsilon_{k\ell}}(s)
(\cK_{\ell} m_k)(s)ds\right) \xi_{\alpha}
\ee
as long as the series in $\alpha$ converges in
the mean square (the inner sum always contains
finitely many non-zero terms).
\end{proposition}

\begin{proof} Part (a) is a classical result of Cameron
and Martin \cite{CM}. Part (b) follows from the
definition of the It\^{o}-Skorokhod integral
in terms of the Wick product \cite{LtSt}.

\end{proof}

Let us now assume that
 equation \eqref{eq:LSp-eqm1}
has a variational solution $u$ and
\bel{eq:sol-chaos}
u(t)=\sum_{\alpha\in \cJ} u_{\alpha}(t)\xi_{\alpha}.
\ee
Substituting this representation into
\eqref{eq:LSp-eqm1} and using \eqref{eq:int-chaos},
we conclude that each $u_{\alpha}$, which is
non-random, satisfies
 \bel{s-system}
\begin{split}
u_{\alpha}(t)&=
u_{0,\alpha} +
\int_0^t\big(\cA u_{\alpha}(s)
+ F_{\alpha}(s)\big)ds\\
&+ \sum_{k,\ell=1}^{\infty}\sqrt{\alpha_{k\ell}}
\int_0^t\big(\cB_{\ell} u_{\alpha-\epsilon_{k\ell}}(s)
+G_{\ell,\alpha-\epsilon_{k\ell}}(s))
(\cK_{\ell} h_k)(s)ds,
\end{split}
\ee
where $u_{0,\alpha}=\bE(u_0\xi_{\alpha})$,
$F_{\alpha}(t)=\bE(F(t)\xi_{\alpha})$,
$G_{\ell,\alpha}(t)=\bE(G_{\ell}(t)\xi_{\alpha})$.
This observation motivates the following definition
of the chaos solution.

\begin{definition}
\label{def:chaos-sol}
(a) The collection of functions
$\{u_{\alpha},\ \alpha\in \cJ\}$ is called a
\und{chaos solution} of equation \eqref{eq:LSp-eqm1}
if every $u_{\alpha}$ is an element of
$L_2((0,T); \bfX)$ and the system of equalities
\eqref{s-system}
holds in $\bfX'$ for all $t\in [0,T]$.
The chaos solution is called \und{square integrable}
if
\begin{equation}
\sup_{0<t<T}
\sum_{\alpha\in \cJ} \|u_{\alpha}(t)\|_{\bH}^2 < \infty.
\end{equation}
(b) The system of equalities \eqref{s-system} is
called the \und{S-system}
corresponding to equation \eqref{eq:LSp-eqm1}.
\end{definition}

\begin{remark}
(a) If $\{u_{\alpha},\ \alpha \in \cJ\}$ is a
square integrable chaos solution, then, for
each $t\in [0,T]$, $u(t)=\sum_{\alpha \in \cJ}
u_{\alpha}(t)\xi_{\alpha}$ is an element of
$L_2(\Omega; \bH)$, but still there is no guarantee
that $u$ is a variational solution.

(b) Uniqueness of the chaos solution implies uniqueness
of the variational solution.
\end{remark}

To establish existence and uniqueness of the chaos
solution, we look at \eqref{eq:LSp-eqm1} as
a system of equations. To solve this system,
 we  make the following assumptions:
\begin{enumerate}
\item[{\bf A1}]
The operator $\cA$ is bounded linear
from $\bfX$ to $\bfX'$ and is \und{strongly parabolic}:
  there exist a positive number $\delta_A$ and a real
number $C_A$ such that, for all $v\in \bfX$,
\bel{eq:parab-c}
 \langle \cA v,v \rangle +\delta_A\|v\|_{\bfX}^2
  \leq C_A\|v\|_{\bH}^2.
 \ee
\item[{\bf A2}]
Each $\mfX_{\ell},\ \ell\geq 1,$ is a
  colored noise on $L_2((0,T))$
with a representation operator $\cK_{\ell}$
and each $\cK_{\ell}$ is a bounded linear operator on
 $L_2((0,T))$ with the operator norm
 $\mathfrak{K}_{\ell}:$
\begin{equation}
\int_0^T|(\cK_{\ell} f)(s)|^2ds
\leq \mathfrak{K}_{\ell}^2 \int_0^T f^2(s)ds,\ \
f\in L_2((0,T)).
\end{equation}
\item[{\bf A3}] The functions
 $u_0, \,F, \,G_{\ell}, \ \ell\geq 1$, are non-random
 and  $\mfX_{\ell},\ \ell\geq 1,$
 are jointly independent.
\item[{\bf A4}]
 Each $\cB_{\ell}$ is a bounded linear
operator on $\bH$ with the operator norm $C_{\ell}$:
\begin{equation}
\|\cB_{\ell} v\|_{\bH}
\leq C_{\ell} \|v\|_{\bH},\ v\in \bH.
\end{equation}
\item[{\bf A5}]
The following regularity conditions hold:
\bel{eq:input}
I_0=\|u_0\|_{\bH}^2+
\frac{2}{\delta_A}\int_0^T \|F(t)\|^2_{\bfX'}dt
+\sum_{\ell=1}^{\infty} \mathfrak{K}_{\ell}^2
\int_0^T\|G_{\ell}(t)\|^2_{\bH}dt < \infty,
\ee
\bel{eq:oper}
C_B=\sum\limits_{\ell=1}^{\infty}
\mathfrak{K}_{\ell}^2\,C_{\ell}^2<\infty.
\ee
\end{enumerate}

 \begin{theorem}
 \label{th:chaos-m1}
Under Assumptions {\bf A1}--{\bf A5}
 equation \eqref{eq:LSp-eqm1}
has a unique chaos solution. The solution is
square integrable and satisfies
\bel{eq:reg-sol}
\sup_{0<t<T}\bE\|u(t)\|_{\bH}^2
\leq C_oe^{(C_A+C_B)T}I_0,
\ee
where $C_A>0$ is from \eqref{eq:parab-c} and
 $1\leq C_o\leq 3.$
 In particular, $C_o=1$ if $F=G_{\ell}=0.$
\end{theorem}

\begin{proof}
Under Assumption {\bf A3} the
 the S-system \eqref{s-system} corresponding to equation
 \eqref{eq:LSp-eqm1} becomes
\bel{s-system1}
\begin{split}
u_{(0)}(t)&=u_0+\int_0^t \cA u_{(0)}(s)ds
+\int_0^t F(s)ds,\ |\alpha|=0;\\
u_{\epsilon_{ij}}(t)&=
\int_0^t \cA u_{\epsilon_{ij}}(s)ds+
\int_0^t \Big(
\cB_{j}u_{(0)}(s)+G_{j}(s)
\Big)(\cK_{j}h_i)(s)ds,\ |\alpha|=1;\\
u_{\alpha}(t)&=
\int_0^t\cA u_{\alpha}(s)ds
+ \sum_{k,\ell=1}^{\infty}\sqrt{\alpha_{k\ell}}
\int_0^t\cB_{\ell} u_{\alpha-\epsilon_{k\ell}}(s)
(\cK_{\ell} h_k)(s)ds, \ |\alpha|>1.
\end{split}
\ee
Indeed, by Assumption {\bf A3}, if $|\alpha|>0$,
then
 $\bE(u_{0}\xi_{\alpha})=0$, $\bE(F(t)\xi_{\alpha})=0$,
  and $\bE(G_{\ell}\xi_{\alpha})=0$.

 The following proposition provides a key estimate
 for the solution of this system and is the main
 step in the proof of Theorem \ref{th:chaos-m1}.

 \begin{proposition}
 \label{th:conv-main}
 Under Assumptions {\bf A1}--{\bf A6},
 for every $0\le t\leq T$ and $k\geq1$,
\bel{s-main1}
\begin{split}
\sum_{\stackrel{\scriptstyle{\alpha\in\cJ}}{|\alpha|=k}}
\|u_{\alpha}(t)\|_{\bH}^2&\leq
 C_oe^{C_AT} \Bigg( \frac{(C_B\,T)^{k}}{k!}\,
\left(\|u_0\|_{\bH}^2+
\frac{2}{\delta_A}\int_0^T\|F(s)\|_{\bfX'}^2ds\right)\\
&+\frac{(C_BT)^{k-1}}{(k-1)!}
\sum_{\ell=1}^{\infty}\mathfrak{K}^2_{\ell}
\int_0^T \|G_{\ell}(s)\|_{\bH}^2ds\Bigg).
\end{split}
\ee
\end{proposition}

\begin{proof}
Given $v_0\in \bH$ and
$\psi\in L_2((0,T); \bfX')$, consider a
deterministic evolution equation
\bel{eq:det0}
v(t)=v_0+\int_0^t(\cA v(s) +\psi(s))ds.
\ee
By definition, $v\in L_2((0,T);\bfX)$ is a solution
of \eqref{eq:det0} if equality \eqref{eq:det0}
holds in $\bfX'$ for every $t\in [0,T]$.

 It is known \cite[Theorem 3.1.4]{Roz}
 that, if the operator $\cA$ is strongly
 parabolic, then
 \begin{itemize}
 \item The operator $\cA$ generates a
 semigroup $\Phi=\Phi_t, \, t\geq 0$, in the space $\bH.$
\item The semigroup $\Phi$ has the following properties:
 \bel{eq:3.2}
 \|\Phi_t v\|_{\bH}^2\leq e^{C_At}\|v\|_{\bH}^2,
 \ v\in \bH;
 \ee
 \bel{eq:3.3}
 \left\|\int_0^t\Phi_{t-s} f(s)ds\right\|_{\bH}^2
 \leq \frac{2}{\delta_A}\,
 e^{C_At}\int_0^t\|f(s)\|_{\bfX'}^2ds,
 \ f\in L_2((0,T);\bfX').
 \ee
\item The solution of  (\ref{eq:det0}) is unique
 and  can be written as
 \bel{eq:3.4}
  v(t)=\Phi_tv_0+\int_0^t\Phi_{t-s}\psi(s)ds.
 \ee
   \end{itemize}

We use these results to study the system of
equations (\ref{s-system1}).
It follows by induction on $|\alpha|$ that
 if $|\alpha|=k$ with the characteristic set
 $\{(i_1, \ell_1), \ldots, (i_k,\ell_k)\},$
 and $\mathcal{P}_k$ is the set of all permutations
 of $\{1,2,\ldots,k\}$,
 then the solution of (\ref{s-system1})
 is unique and is given by
 \begin{equation}
\begin{split}
\label{eq:ind_sg}
u_{\alpha}(t)=
&\frac{1}{\sqrt{\alpha!}}\sum_{\sigma \in \cP_k}
\int_0^t\int_0^{s_{k}}\ldots \int_0^{s_2}
\Phi_{t-s_k}\cB_{\ell_{\sigma(k)}}\cdots
\Phi_{s_2-s_1}\Big(\cB_{\ell_{\sigma(1)}}u_{(0)}(s_1)\\
&+G_{\ell_{\sigma(1)}}(s_1)\Big)
(\cK_{\ell_{\sigma(k)}} h_{i_{\sigma(k)}})(s_k) \cdots
(\cK_{\ell_{\sigma(1)}} h_{i_{\sigma(1)}})(s_1)
ds_1\ldots ds_k.
\end{split}
\end{equation}
We then re-write (\ref{eq:ind_sg}) as
\bel{eq:ind_sg1}
u_{\alpha}(t)=
\int_{[0,T]^k}H(t,\ell^{(k)};s^{(k)})
\overline{h}_{\alpha}(s^{(k)})ds_1 \ldots ds_k,
\ee
where
\begin{equation}
\begin{split}
H(t,\ell^{(k)};s^{(k)})
&=\frac{1}{\sqrt{k!}}\sum_{\sigma\in \cP^k}
\Phi_{t-s_{\sigma(n)}}\cB_{\ell_n}\cdots
\Phi_{s_{\sigma(2)}-s_{\sigma(1)}}
\Big( \cB_{\ell_1}u_{(0)}(s_{\sigma(1)})\\
&+G_{\ell_1}(s_{\sigma(1)})\Big)
\chi_{s_{\sigma(2)}}(s_{\sigma(1)})
\cdots\chi_t(s_{\sigma(k)}),
\end{split}
\end{equation}
and
\begin{equation}
\overline{h}_{\alpha}(s^{(k)})
=\frac{1}{\sqrt{\alpha!\,k!}}\sum_{\sigma \in\cP_k}
 (\cK_{\ell_1} h_{i_1})(s_{\sigma(1)})
\cdots (\cK_{\ell_k} h_{i_k})(s_{\sigma(k)}).
\end{equation}
From (\ref{eq:ind_sg1}) and the definition of the
function $H$, we conclude that
\begin{equation}
\label{mult-int}
\begin{split}
\sum_{\stackrel{\scriptstyle{\alpha\in\cJ}}{|\alpha|=k}}
&\|u_{\alpha}(t)\|^2_{\bH}
\leq \sum_{\ell_1,\ldots,\ell_k=1}^{\infty}
\left(\prod_{j=1}^k \mathfrak{K}_{\ell_j}^2\right)
 \int_0^t\int_0^{s_{k}}\ldots \int_0^{s_2}
 \\ & \left\|\Phi_{t-s_k}\cB_{\ell_{k}}\cdots
\Phi_{s_2-s_1}\Big(\cB_{\ell_{1}}u_{(0)}(s_1)\right.
\left.+G_{\ell_1}(s_1)\Big)\right\|_{\bH}^2
ds_1\ldots ds_k.
\end{split}
\end{equation}
By \eqref{s-system1} and \eqref{eq:3.4},
\begin{equation}
u_{(0)}(t)=\Phi_tu_0+\int_0^t\Phi_{t-s}F(s)ds,
\end{equation}
and then the
properties \eqref{eq:3.2} and \eqref{eq:3.3} of the
semigroup $\Phi$ imply
$$
\|\cB_{\ell}u_{(0)}(t)+G_{\ell}(t)\|_{\bH}^2
\leq C_o \Big( C_{\ell}^2e^{C_At} \|u_0\|_{\bH}^2
+\frac{2}{\delta_A}\,C_{\ell}^2 e^{C_At}
\int_0^t \|F(s)\|^2_{\bfX'}ds
+\|G_{\ell}(t)\|_{\bH}^2 \Big),
$$
where $C_o=3$, as the inequality
$(a+b+c)^2\leq 3(a^2+b^2+c^2)$ suggests.
On the other hand,
$C_o=1$ if $F=G_{\ell}=0$.
Applying \eqref{eq:3.2} repeatedly, we find
\begin{equation}
\begin{split}
&\|\Phi_{t-s_k}\cB_{\ell_{k}}\cdots
\Phi_{s_2-s_1}\big(\cB_{\ell_{1}}u_{(0)}(s_1)
+G_{\ell_1}(s_1)\big)\|_{\bH}^2\\
&\leq C_{\ell_k}^2e^{C_A(t-s_k)}
\|\Phi_{s_k-s_{k-1}}\cB_{\ell_{k-1}}\cdots
\Phi_{s_2-s_1}\big(\cB_{\ell_{1}}u_{(0)}(s_1)
+G_{\ell_1}(s_1)\big)\|_{\bH}^2
\leq \ldots\\
& \leq
C_o\left(\prod_{j=2}^kC_{\ell_j}^2\right)e^{C_At}
\left( C_{\ell_1}^2\|u_{0}\|_{\bH}^2
+\frac{2}{\delta_A}\,C_{\ell_1}^2
\int_0^T\|F(s)\|_{\bfX'}^2ds+
\|G_{\ell_1}(s_1)\|_{\bH}^2\right).
\end{split}
\end{equation}
Inequality (\ref{s-main1}) now follows from
\eqref{mult-int}.

Proposition \ref{th:conv-main} is proved.
\end{proof}

To complete the proof of Theorem \ref{th:chaos-m1}
it remains to note that uniqueness of the
chaos solution is equivalent to the uniqueness
of the solution of \eqref{s-system1} and is
guaranteed by the
strong parabolicity of the operator $\cA$.
Then
\begin{equation}
\sum_{\alpha \in \cJ} \|u_{\alpha}(t)\|_{\bH}^2
=\|u_{(0)}(t)\|_{\bH}^2
+\sum_{k\geq 1}
\sum_{\stackrel{\scriptstyle{\alpha\in\cJ}}{|\alpha|=k}}
\|u_{\alpha}(t)\|_{\bH}^2,
\end{equation}
and \eqref{eq:reg-sol} follows from \eqref{s-main1}.
\end{proof}

As a first step toward studying the approximation
of the chaos solution, we get
\begin{corollary}
\label{cor-appr0}
Let $F=G_{\ell}=0$ and, for $N\geq 1$, define
\bel{appr1}
u_N(t)=
\sum_{\stackrel{\scriptstyle{\alpha\in\cJ}}
{|\alpha|\leq N}}
u_{\alpha}(t)\xi_{\alpha}.
\ee
Then
\bel{appr2}
\sup_{0<t<T} \bE \|u(t)-u_N(t)\|_{\bH}^2
\leq \frac{(C_BT)^{N+1}}{(N+1)!} e^{(C_A+C_B)T}
 \|u_0\|_{\bH}^2.
\ee
\end{corollary}

\begin{proof}
We have
\begin{equation}
\bE \|u(t)-u_N(t)\|_{\bH}^2=
\sum_{\stackrel{\scriptstyle{\alpha\in\cJ}}
{|\alpha|\geq N+1}}
\|u_{\alpha}(t)\|_{\bH}^2,
\end{equation}
because
\bel{eq:xi-orth}
\bE(\xi_{\alpha}\xi_{\beta})=
\begin{cases}
1& {\rm if} \ \alpha=\beta,\\
0& {\rm otherwise}.
\end{cases}
\ee
 Then \eqref{appr2} follows
from \eqref{s-main1}.
\end{proof}

\begin{remark} The conclusions of the
theorem are valid even if the operators
$\cA$, $\cB_{\ell}$ depend on time in a sufficiently
regular way, as long
as Assumptions {\bf A1} and {\bf A3} hold uniformly
in $t\in [0,T]$.
\end{remark}

There are at least two  open problems
related to the chaos solution of
equation \eqref{eq:LSp-eqm1}:
\begin{enumerate}
\item To find out whether the solution belongs to
\begin{equation*}
L_2(\Omega; L_2(0,T);X))\bigcap
L_2(\Omega; C((0,T);H));
\end{equation*}
  this is true when
every $\mfX_{\ell}$ is a white noise over
$L_2((0,T))$ \cite[Theorem 3.8]{LR1}.
\item To establish existence of the
solution when the operators $\cB_{\ell}$ are unbounded
on $\bH$.
\end{enumerate}

\section{Approximation of the Chaos Solution}
\setcounter{equation}{0}
\label{Sec:appr}

\subsection{One-step Approximation}

\subsubsection{Motivation}

If $\bA$ is a Hilbert space,
 $f\in \bA$, and $\{m_k,\  k\geq 1\}$ is an
orthonormal basis in $\bA$, then
$$
\sum_{k\geq 1}|(f,m_k)_{\bA}|^2<\infty,
$$
but nothing can be said about the rate of
this convergence, that is, about the rate at which
$$
\sum_{k=n}^{\infty}|(f,m_k)_{\bA}|^2
$$
tends to zero as $n\to \infty:$
taking $\bA=L_2((0,T))$
with a trigonometric basis, one can construct
a function for which this convergence will be
arbitrarily slow.

In the study of the chaos solution, we are facing a
similar problem.
In fact, the underlying Hilbert space is the
 space of square integrable random processes,
 the study of the rate of convergence for
the ``natural'' approximations of the chaos solution is
reduced to the analysis of certain Fourier series
in $L_2((0,T))$.

As an illustration, consider the
following equation:
\bel{spde-cp333}
u(t,x)=u_0+\int_0^tu_{xx}(s,x)ds
+\mfX(\chi_t\,h(\cdot,x)u(\cdot,x)),\
t\geq 0,\ x\in \bR,
 \ee
where $\mfX$ is a colored noise on $L_2((0,T))$
with representation $(\cK,W)$.

 If $u_0\in L_2(\bR)$
is non-random, and
$h=h(t,x)$ is a bounded non-random function,
then Theorem \ref{th:chaos-m1} on
page \pageref{th:chaos-m1}
implies that
(\ref{spde-cp333}) has a unique square-integrable
chaos solution  $u(t,x)=\sum\limits_{\alpha\in \cI}
u_{\alpha}(t,x)\xi_{\alpha}$, where
\bel{s-aux333}
\begin{split}
\frac{\partial u_{(0)}}{\partial t}
&=\frac{\partial^2 u_{(0)}}{\partial x^2},\
 u_{(0)}(0,x)=u_0(x),\\
\frac{\partial u_{\alpha}}{\partial t}
&=\frac{\partial^2 u_{\alpha}}{\partial x^2}+
\sum_{k=1}^{\infty}\sqrt{\alpha_k}\,
hu_{\alpha-\epsilon_k}\cK m_k,\
u_{\alpha}(0,x)=0, \ |\alpha|>0,
\end{split}
\ee
and $\{m_k,\, k\geq 1\}$ is an orthonormal basis
in $L_2((0,T))$; with only one colored noise driving
the equation, every multi-index has the form
$\alpha=\{\alpha_1,\alpha_2, \ldots\}$.

Denoting the heat semigroup by $\Phi_t$, we find
 \bel{s-aux3332}
u_{(0)}(t)=\Phi_tu_0(x)
\ee
and
\bel{s-aux3333}
u_{\epsilon_k}(t,x)
=\int_0^t\Phi_{t-s}h \Phi_su_{0}(s,x)(\cK m_k)(s)ds.
\ee
Let us define an approximation $u_1^n(t,x)$
of $u(t,x)$ by
\bel{ex-appr0}
u_1^n(t,x)=u_{(0)}(t,x)
+\sum_{k=1}^nu_{\epsilon_k}(t,x)\xi_k.
\ee
What can we say about the quality of this
 approximation? For example, can
we find a bound on
$\sup_{0<t<T}\bE\|u-u_1^n\|_{L_2(\bR)}^2(t)$ in terms of
$n$ and $T$?

Since
$u(t,x)=\sum\limits_{\alpha\in \cI}
u_{\alpha}(t,x)\xi_{\alpha}$
and \eqref{eq:xi-orth} holds, we have
\bel{ex-appr2}
\bE\|u-u_1^n\|_{L_2(\bR)}^2(t)=
\sum_{\alpha\in \cI,\ |\alpha|>1}
\|u_{\alpha}\|_{L_2(\bR)}^2(t)
+\sum_{k=n+1}^{\infty}
\|u_{\epsilon_k}\|_{L_2(\bR)}^2(t).
\ee
 Using the properties of the heat semigroup on
 $\bR$ and Corollary \ref{cor-appr0} on page
  \pageref{cor-appr0}, we find
\bel{ex-appr1}
\sup_{0<t<T}\sum_{\alpha\in \cI,
\ |\alpha|>1} \|u_{\alpha}\|_{L_2(\bR)}^2(t)
\leq (C_BT)^2e^{C_BT},
\ee
where  $C_B=\|\cK\|\,\sup\limits_{t,x}|h(t,x)|^2$.

As a result, to find the quality of
the approximation, we need to
find the  rate of convergence of the series
$
\sum_{k=1}^{\infty} \|u_{\epsilon_k}\|_{L_2(\bR)}^2(t).
$
This rate of convergence is determined by the
rate of decay, as
$k\to \infty$, of $\|u_{\epsilon_k}\|_{L_2(\bR)}^2(t)$,
 and, as equality (\ref{s-aux3333}) suggests,
 one way to determine this rate is to
integrate by parts. Accordingly, setting
\bel{M000}
\wt{M}_k(t)=\int_0^t (\cK m_k)(s)ds,
\ee
and using the properties of the
heat semi-group, we get
\bel{conv-fc000}
\begin{split}
u_{\epsilon_k}(t,x)&=
\Phi_{t-s}(h(s,\cdot)u_{(0)}(s,\cdot))(x)
\wt{M}_k(s)\Big|_{s=0}^{s=t}\\
&-\int_0^t\Big(\Phi_{t-s}h(\Phi_s u_{0})_{xx}
-\big(\Phi_{t-s}h\Phi_su_0\big)_{xx}
\Big)\wt{M}_k(s)ds.
\end{split}
\end{equation}

Note that (\ref{M000}) and the Cauchy-Schwartz
inequality imply
\begin{equation}
|\wt{M}_k(t)|\leq C\sqrt{t},
\end{equation}
and so $\wt{M}_k(0)=0$. Still, to advance our study
of the rate of convergence
any further, we need
\begin{enumerate}
\item additional regularity of $u_0$ and $h$;
\item a rather detailed information about the
functions $\wt{M}_k$.
\end{enumerate}

If we indeed assume all the necessary regularity
of $u_0$ and $h$, then (\ref{conv-fc000}),
together with the Cauchy-Schwartz inequality, implies
 \bel{conv-fc001}
 \|u_{\epsilon_k}\|_{L_2(\bR)}^2(t)
 \leq C_1\|u_0\|_{L_2(\bR)}^2|\wt{M}_k(t)|^2+
C_2
\left\|\frac{\partial^2u_{0}}{\partial x^2}
\right\|_{L_2(\bR)}^2\ t\int_0^t|\wt{M}_k(s)|^2ds.
\ee
To continue, assume that
\begin{equation}\label{basis}
 m_1(s)\!=\!\frac 1{\sqrt{T}};\  m_k(t)\!=
 \!\sqrt{\frac{2}{T} }
 \cos \left( \frac{\pi (k-1) t}{T} \right), \, k>1; \
0\leq t \leq T.
\end{equation}
With this choice of the basis, we use
\eqref{M000} to find that, for $k>1$,
\begin{enumerate}
\item If $\mfX$ is white noise ($\cK=I$), then
\bel{M1-1}
\wt{M_k}(t)
=\frac{\sqrt{2T}}{\pi(k-1)}
\sin\left(\frac{\pi (k-1) t}{T} \right);
\ee
\item If $\mfX$ is the
Ornstein-Uhlenbeck noise with parameter $b$, then
\bel{M1-2}
\begin{split}
\wt{M}_k(t)&=\frac{\sqrt{2T^3}}{b^2T^2
+(k-1)^2\pi^2}
\Bigg(b \cos \left( \frac{\pi (k-1) t}{T} \right)
-be^{-bt}\\
&+\frac{(k-1)\pi}{T}\sin \left( \frac{\pi (k-1) t}{T}
 \right)\Bigg);
\end{split}
\ee
\item If $\mfX$ the $H$-{fractional white noise}
 and $1/2<H<1$, then
\bel{M1-3}
|\wt{M}_k(t)|
\leq \frac{C(H)t^{2H-1}T^{1-H}}{k^{\frac{3}{2}-H}}
\ee
for some number $C(H)$ depending only on $H$.
\end{enumerate}
Relation (\ref{conv-fc001}) suggests
that the rate of convergence
will be quite different for different $\mfX$,
and we are essentially forced
 to make the following  assumptions about
the functions $\wt{M}_k$:
\begin{align}
\sup_{0<t<T} |\wt{M}_k(t)|^2
\leq \wt{C}\ \frac{T^{\delta}}{k^{\gamma}} \quad
{\rm for \ } \delta>0,\ \gamma>1,
\  \wt{C}>0\label{m-ass-1},\\
|\wt{M}_k(T)|^2
\leq \wt{C}\ \frac{T^{\delta_1}}{k^{\gamma_1}} \quad
{\rm for \ } \delta_1>0,\
\gamma_1>1, \ \wt{C}>0 \label{m-ass-2}.
\end{align}
In both (\ref{m-ass-1}) and (\ref{m-ass-2}),
the number $\wt{C}$ should  not
depend on $T$ or $k$.

It is enough to have (\ref{m-ass-1}) and (\ref{m-ass-2})
 for {\em some} orthonormal basis $\{m_k,\ k\geq 1\}$
in $L_2((0,T))$, but for now the
cosine basis (\ref{basis}) is
the only example when these assumptions can be
verified. Below, we summarize
the results for the white noise $W$, fractional
white noise
$W^H$, $1/2<H<1$, and the Ornstein-Uhlenbeck noise
 $U_b$, $b>0$,
when the cosine basis (\ref{basis}) is used.\\

\begin{table}
\caption{}
\label{tab:rate}
\begin{tabular}{ccccc}
\hline\hline
& & & &\\
    &                 $\delta$ & $\gamma$ &  $\delta_1$ & $\gamma_1$ \\
& & & & \\ \hline
& & & & \\
$W$                   & $1$   & $2$       & Any & Any   \\
&&&& \\ \hline
& & & & \\
$W^H$  & $\ 2H\ $  &\  $3-2H\ $   & $\ \ 2H\ \ $  & $\ 3-2H\ $ \\
& & & & \\ \hline
& & & & \\
$U_b$            & $1$  & $2$      & $3$   & $4$ \\
& & & & \\ \hline \hline
\end{tabular}
\end{table}

\begin{remark}
(a) In the case of $W$, with $m_k$
 as in (\ref{basis}),
we have $\wt{M_k}(T)=0$ for all $ k\geq 2$, and then
indeed any choice of $\delta_1$, $\gamma_1$ will work in
\eqref{m-ass-2}.

(b) Inequality \eqref{conv-fc001} shows that assumptions
(\ref{m-ass-1}) and (\ref{m-ass-2})
are close to necessary for the analysis
of convergence of the
chaos expansion.

(c) As equality
(\ref{conv-fc000}) suggests, no further
 integration by parts will, in
general, improve the rate of convergence.
\end{remark}

Under assumption (\ref{m-ass-1}),  we conclude from
(\ref{conv-fc001}) that, for $n>1$,
\bel{conv-fc005}
\sup_{0<t<T}\sum_{k=n}^{\infty}
\|u_{\epsilon_k}\|_{L_2(\bR)}^2(t)
\leq \wt{C}_1 \|u_0\|_{L_2(\bR)}^2
 \frac{T^{\delta}}{n^{\gamma-1}}+
\wt{C}_2 \left\|\frac{\partial^2u_{0}}{\partial x^2}
\right\|_{L_2(\bR)}^2
 \frac{T^{\delta+2}}{n^{\gamma-1}}.
\ee
With assumption (\ref{m-ass-2}), we also get
\bel{conv-fc008}
\sum_{k=n}^{\infty}
\|u_{\epsilon_k}\|_{L_2(\bR)}^2(T)
\leq \wt{C}_1 \|u_0\|_{L_2(\bR)}^2
 \frac{T^{\delta_1}}{n^{\gamma_1-1}}
+\wt{C}_2 \left\|\frac{\partial^2u_{0}}{\partial x^2}
\right\|_{L_2(\bR)}^2
\frac{T^{\delta+2}}{n^{\gamma-1}}.
\ee

Inequality  (\ref{conv-fc005}) establishes
an approximation error bound
uniformly over the time
interval $(0,T)$, while (\ref{conv-fc008})
gives the bound only at the end
point. If $T$ is small, and if we can take
 $\delta_1>\delta$, $\gamma_1>\gamma$,
which is the case for $\mfX=\mfB$ and
$\mfX=U_b$, then (\ref{conv-fc008}) provides
 a better error bound
than (\ref{conv-fc005}) and is
more suitable for analyzing a step-by-step
approximation.

By combining (\ref{ex-appr1}) with either
 (\ref{conv-fc005}) or
(\ref{conv-fc008}), we will get the overall bound on the
approximation error; for white, fractional, or
Ornstein-Uhlenbeck noise, we also use
Table \ref{tab:rate}. For example,
when $\mfX=W^H$, $1/2<H<1$, we
have $\|\cK\|^2\leq C(H)T^{2H-1}$
(see Example \ref{Exe-fBm}, page \pageref{Exe-fBm})
 and therefore, for $T\leq 1$,
\bel{ex-appr-fbm}
\sup_{0<t<T}\bE\|u-u_1^n\|_{L_2(\bR)}^2(t)\leq
C^{*} \left(T^{4H}+\frac{T^{2H}}{n^{2-2H}}\right).
 \ee
 In the next
section, we extend this result to more general
 equations and more
general approximations.

\subsubsection{Truncation of the S-system}
\label{sec-tr-s}

Consider the following evolution equation:
 \bel{tr-evol1}
u(t)=u_0+\int_0^t \cA u(s) ds
 + \sum_{\ell=1}^{\infty}
\mfX_{\ell}(\chi_t\cB_{\ell} u),
 \ee
 and assume that this
equation has a unique square-integrable chaos
solution $u$ in a normal
triple $(\bfX,\bH,\bfX')$ of Hilbert spaces.
 We also assume that $u$
has chaos expansion
\bel{tr-sol1} u(t)=\sum_{\alpha\in \cJ}
u_{\alpha}(t)\xi_{\alpha}.
 \ee
 As before, we assume that $u_0$ is deterministic and
 every colored noise $\mfX_{\ell}$ has a
 representation $(\cK_{\ell}, W_{\ell})$.

The fist step
is to find a general method of constructing an
approximation of $u$ given
the expansion (\ref{tr-sol1}).
A natural approximation is
\bel{tr-sol2}
\ol{u}(t)=\sum_{\alpha\in \ol{\cJ}}
u_{\alpha}(t)\xi_{\alpha},
\ee
where $\ol{\cJ}$ is a
{\em finite} subset of $\cJ$. To control the
size of this finite set, we use
three characteristics of a
multi-index:
$$
|\alpha|=\sum_{k,\ell}\alpha_{k\ell},
\ \varpi(\alpha)= \max\{k:\alpha_{k\ell}>0\},\
d(\alpha)=\max\{\ell:\alpha_{k\ell}>0\}.
$$
For example, if
$$
\alpha=
\left(
\begin{array}{lllllllll}
1 & 0 & 1 & 0 & 0 & 0 & 3 & 0 & 0 \cdots\\
0 & 0 & 0 & 1 & 0 & 0 & 0 & 2 & 0 \cdots\\
1 & 0 & 0  & 0& 0 & 0 & 0 & 0& 0 \cdots \\
0 & 0 & 0  & 0& 0 & 0 & 0 & 0 & 0 \cdots \\
\vdots &  \vdots & \vdots &  \vdots & \vdots &  \vdots & \vdots &  \vdots & \vdots
\end{array}
\right),
$$
then  $|\alpha|=1+1+3+1+2+1=9$, $\varpi(\alpha)=8$,
 $d(\alpha)=3$.
We call $|\alpha|$ the \und{length of the multi-index},
$\varpi(\alpha)$ the \und{order of the multi-index}, and
$d(\alpha)$ the \und{dimension of the multi-index}.
 Then the set
$$
\cJ^{n,r}_N=\{\alpha\in \cJ: |\alpha|\leq N,\
\varpi(\alpha)\leq n,\ d(\alpha)\leq r\}
$$
is finite, with no more than $(nr)^N$ elements.
 Note that the sets
$\cJ_N=\{\alpha\in \cJ: |\alpha|\leq N\}$
is always infinite, and the set
$\cJ^n_N=\{\alpha\in \cJ:
|\alpha|\leq N, \varpi(\alpha)\leq n\}
$ is infinite if and only if there are
infinitely many noises in the equation.
Accordingly, we define three approximations of $u$:
\bel{appr-gen}
u_N(t)=\sum_{\alpha\in \cJ_N}u_{\alpha}(t)\xi_{\alpha},
    \ u_N^n(t)
    =\sum_{\alpha\in \cJ_N^n}u_{\alpha}(t)\xi_{\alpha},\
 u_N^{n,r}(t)
 =\sum_{\alpha\in \cJ_N^{n,r}}u_{\alpha}(t)\xi_{\alpha}
 \ee
Of the three, only $u_N^{n,r}(t)$ is {\em computable},
 being a sum of finitely many terms.
  Consequently, our goal is to find a bound
 on
 $\sup\limits_{0<t<T}\bE\|u(t)-u_N^{n,r}(t)\|_{\bH}$.
  Recall
 (see (\ref{s-system}) on page \pageref{s-system})
 that the coefficients $u_{\alpha}$ satisfy the S-system
 \bel{err-s-syst}
 \begin{split}
&u_{(0)}(t)=u_0+\int_0^t \cA u_{(0)}(s)ds,\ \
|\alpha|=0; \\
& u_{\alpha}(t)=
\int_0^t\cA u_{\alpha}(s) ds
+
\sum_{k,\ell=1}^{\infty}\sqrt{\alpha_{k\ell}}
\int_0^t\cB_{\ell} u_{\alpha-\epsilon_{k\ell}}(s)\,
(\cK m_k)(s)ds,\ |\alpha|>0.
\end{split}
 \ee
  By orthogonality
 of $\xi_{\alpha}$ for different $\alpha,$
 we have
 \begin{equation}
 \label{orth-000}
 \begin{split}
\bE\|u(t)-u_N^{n,r}(t)\|_{\bH}^2
&=\sum_{\alpha\notin \cJ_N^{n,r}}
\bE\|u_{\alpha}(t)\|_{\bH}^2 =
\sum_{\alpha\in \cJ\backslash \cJ_N}
\bE\|u_{\alpha}(t)\|_{\bH}^2\\
&+
\sum_{\alpha\in \cJ_N\backslash\cJ_N^n}
\bE\|u_{\alpha}(t)\|_{\bH}^2+
\sum_{\alpha\in \cJ_N^n\backslash\cJ_N^{n,r}}
\bE\|u_{\alpha}(t)\|_{\bH}^2,
\end{split}
\end{equation}
where $\backslash$ denotes the difference of two sets.
In other words,
we have an analogue of the Pithagorean theorem:
\bel{orth-001}
\begin{split}
\bE\|u(t)-u_N^{n,r}(t)\|_{\bH}^2&=
\bE\|u(t)-u_N(t)\|_{\bH}^2+
\bE\|u_N(t)-u_N^{n}(t)\|_{\bH}^2\\
&+ \bE \|u_N^n(t)-u_N^{n,r}(t)\|_{\bH}^2,
\end{split}
\ee
and Corollary \ref{cor-appr0}, page
\pageref{cor-appr0}, provides an estimate for
$\bE\|u(t)-u_N(t)\|_{\bH}^2$.

As we saw in the previous section, to estimate
$\bE\|u_N(t)-u_N^{n}(t)\|_{\bH}^2$, we need to assume
(\ref{m-ass-1}) and (\ref{m-ass-2})
together with  additional regularity of the
initial condition $u_0$
and the operators $\cA$, $\cB_{\ell}$.
 To formulate this regularity
we need some additional constructions.

Let $\bH^{r},\ r\in \bR$ be a scale
 of Hilbert spaces,
or a \und{Hilbert scale}
\cite[Section VI.1.10]{KPS} with the
property that $\bH^0=\bH$, $\bH^1=\bfX$,
$\bH^{-1}=\bfX'$. A typical example
of such a scale is the collection of the Sobolev spaces
\bel{sobolev}
\bH^{r}(\bR^d)=\left\{f:
\int_{\bR^d} |\hat{f}(y)|^2(1+|y|^2)^{r}dy<
\infty\right\},
\ee
where $\hat{f}$ is the Fourier transform of $f$.

To generalize the computations
that lead to (\ref{conv-fc001}),
 we make the following assumptions:
\begin{align}
&\bE\|u_0\|_{\bH^2}^2< \infty, \label{ad-reg1}\\
&\|\cA v\|_{\bH}^2 \leq C_{02}\|v\|_{\bH^2}^2, \
\|\Phi_t v\|_{\bH^{j}}^2
\leq e^{C_A t}\|v\|_{\bH^j}^2,\ j=0,2,
\label{ad-reg2}\\
& \cB_{\ell}\cB_{n}
=\cB_n\cB_{\ell} \ {\rm for \ all\ } \ell, \, n,
\label{ad-reg21}\\
&\|\cB_{\ell} v\|_{\bH^2}^2
\leq C_{1,\ell}^2\|v\|_{\bH^2}^2, \ v\in \bH^2,
{\rm\  and\ } \sum_{\ell=1}^{\infty}
 C_{1,\ell}^2 \mathfrak{K}_{\ell}^2=C_{1,B}<\infty,
\label{ad-reg3}
\end{align}
where $\Phi_t$ is the semi-group generated by $\cA$.

\begin{theorem}
\label{th:err-b-2}
Assume that
\begin{itemize}
\item {\bf A1}--{\bf A5} hold (see page
\pageref{eq:parab-c});
\item $F(t)=0$ and $G_{\ell}=0$;
\item (\ref{m-ass-1}), (\ref{m-ass-2}) hold
for all $\mfX_{\ell}$ so that  the numbers
$\delta, \delta_1, \gamma, \gamma_1,$ do not
depend on $\ell$;
\item (\ref{ad-reg1})--(\ref{ad-reg3}) hold.
\end{itemize}
Then
\bel{tr-Nn-1}
\begin{split}
\sup_{0<t<T}\bE\|u_N(t)-u_N^{n}(t)\|_{\bH}^2
 &\leq C_{1,B}
e^{(C_A+\ol{C}_B)T}\Bigg(
C_B\frac{T^{\delta}}{n^{\gamma-1}}\bE\|u_0\|_{\bH}^2\\
&+ C_{02}C_{1,B}\frac{T^{\delta+2}}{n^{\gamma-1}}
\bE\|u_0\|_{\bH^2}^2\Bigg)
\end{split}
\ee
and
\bel{tr-Nn-101}
\begin{split}
\bE\|u_N(T)-u_N^{n}(T)\|_{\bH}^2 &\leq C_{1,B}
e^{(C_A+\ol{C}_B)T}\Bigg(
C_B\frac{T^{\delta_1}}{n^{\gamma_1-1}}
\bE\|u_0\|_{\bH}^2\\
&+
C_{02}C_{1,B}\frac{T^{\delta+2}}{n^{\gamma-1}}
\bE\|u_0\|_{\bH^2}^2\Bigg),
\end{split}
\ee
where $\ol{C}_B=\max(C_{B}, C_{1,B})$.
\end{theorem}

\begin{proof}
 The argument is based on integration by parts in
the representation of $u_{\alpha}(t)$
(see (\ref{eq:ind_sg1}) on page
\pageref{eq:ind_sg1}); while the idea and the
end result are essentially  identical to
(\ref{conv-fc000}), the computations are rather long.
An interested reader can recover these computations
following \cite{LMR}, where each $\mfX$ is white noise.
\end{proof}

\begin{remark}
 Similar to \cite{LMR}, a  bound on
 $\bE\|u_N(t)-u_N^{n}(t)\|_{\bH}^2$ can  be derived
 without condition \eqref{ad-reg21}, that is,
 if the operators $\cB_{\ell}$ do not commute.

 \end{remark}

 Finally, we derive a bound on
 $\bE\|u_N^n(t)-u_N^{n,r}(t)\|_{\bH}^2.$

 \begin{theorem}
 \label{th:err-b-3}
Assume that
\begin{itemize}
\item {\bf A1}--{\bf A5} hold;
\item $F(t)=0$ and $G_{\ell}=0$;
\end{itemize}
Define the sequence $\varepsilon=\varepsilon(r),\
r=1,2,\ldots$ by
 \bel{tr-dim1}
 \sum_{\ell=r+1}^{\infty}
 \mathfrak{K}_{\ell}^2C_{\ell}^2 = \varepsilon(r).
 \ee
 Then
\bel{tr-dim2}
\sup_{0<t<T} \bE \|u_N^n(t)-u_N^{n,r}(t)\|_{\bH}^2
\leq \varepsilon(r)Te^{(C_A+C_B)T}\bE\|u_0\|_{\bH}^2.
\ee
\end{theorem}

\begin{proof} We have by (\ref{mult-int})
 on page \pageref{mult-int}
\begin{equation*}
\begin{split}
\bE \|u_N^n(t)-u_N^{n,r}(t)\|_{\bH}^2 &\leq
\sum_{k=1}^N\sum_{j=1}^k\sum_{\ell_j=r+1}^{\infty}
\sum\nolimits_{\ell_1,\ldots,\ell_k\geq 1}^{(j)}
\int_0^t\int_0^{s_{k}}\ldots \int_0^{s_2}\\
&\|\Phi_{t-s_k}\cB_{\ell_{k}}\cdots
\Phi_{s_2-s_1}\cB_{\ell_{1}}u_{(0)}(s_1)\|_{\bH}^2
ds_1\ldots ds_k,
\end{split}
\end{equation*}
where the summation
$\sum^{(j)}_{\ell_1,\ldots,\ell_k\geq 1}$
omits the index $\ell_j$.
Using assumptions the theorem,
we conclude that
$$
\bE \|u_N^n(t)-u_N^{n,r}(t)\|_{\bH}^2
\leq e^{C_At}
\left(\sum_{\ell=r+1}^{\infty}
\mathfrak{K}_{\ell}^2 C_{\ell}^2
\right)
\left(\sum_{k=1}^{\infty} k\frac{C_B^{k-1}t^k}{k!}
\right)
\bE\|u_0\|^2_{\bH},
$$
which implies (\ref{tr-dim2}).
\end{proof}

\begin{remark}
If there are finitely many $\mfX_{\ell}$,
 and all of them
are included in the approximation,
 then  $\varepsilon(r)=0$ and
 $u_N^n(t)=u_N^{n,r}(t).$
 \end{remark}

Combining the results of Corollary \ref{cor-appr0}
and Theorems \ref{th:err-b-2} and
\ref{th:err-b-3}, we get the overall error bound:
\bel{err-gen-overall}
\begin{split}
\sup_{0<t<T}\bE\|u(t)&-u_N^{n,r}(t)\|_{\bH}^2 \leq
C(T)\Bigg(\frac{(TC_B)^{N+1}}{(N+1)!}\bE\|u_0\|^2_{\bH}\\
&+\frac{T^{\delta}}{n^{\gamma-1}}\bE\|u_0\|^2_{\bH}
+\frac{T^{\delta+2}}{n^{\gamma-1}}\bE\|u_0\|^2_{\bH^2}
+T\varepsilon(r)\bE\|u_0\|^2_{\bH}\Bigg),
\end{split}
\ee
where $\lim\limits_{T\to 0} C(T)>0$.
A similar bound holds for
 $\bE\|u(T)-u_N^{n,r}(T)\|_{\bH}^2$:
 \bel{err-gen-overall1}
\begin{split}
\bE\|u(T)&-u_N^{n,r}(T)\|_{\bH}^2 \leq
C(T)\Bigg(\frac{(TC_B)^{N+1}}{(N+1)!}\bE\|u_0\|^2_{\bH}\\
&+\frac{T^{\delta_1}}{n^{\gamma_1-1}}\bE\|u_0\|^2_{\bH}
+\frac{T^{\delta+2}}{n^{\gamma-1}}\bE\|u_0\|^2_{\bH^2}
+T\varepsilon(r)\bE\|u_0\|^2_{\bH}\Bigg),
\end{split}
\ee

 \begin{example}
Consider the equation
$$
u(t)=u_0+\int_0^tu_{xx}(s)ds+\mfX(\chi_tu),\
0\leq t\leq T,\ x\in \bR.
$$
With only one noise driving the equation, we
have $u_N^{n,r}=u_N^{n}$. Also,
$\bH^{r}=\bH^{r}(\bR)$
is the Sobolev space \eqref{sobolev} and
$\bH=L_2(\bR)$.

 (a) If $\mfX$ is an
 {Ornstein-Uhlenbeck noise}  with
 parameter $b$, then
$C_B=(1+\sqrt{bT})^2,$
 $\delta=1$, $\gamma=2$ (see Table \ref{tab:rate}).
  Inequality  (\ref{err-gen-overall})
 becomes
\bel{err-OU-os}
\begin{split}
\sup_{0<t<T}
\bE\|u(t)-u_N^{n,r}(t)\|_{\bH}^2 &\leq
C_b(T)\Bigg(
\frac{(1+\sqrt{bT})^{2N+2}T^{N}}{(N+1)!}
\bE\|u_0\|^2_{\bH}\\
&+
\frac{T }{n }\bE\|u_0\|^2_{\bH}+
\frac{T^3}{n}\bE\|u_0\|^2_{\bH^2}\Bigg),
\end{split}
\ee
where $\limsup\limits_{T\to 0} C_b(T)>0$.

(b) If $\mfX$ is an
$H$-fractional white noise with
 $H\in (1/2,1)$, then
 $C_B=C_1(H)T^{2H-1},$
 where
 \begin{equation}
 \label{eq:C1H}
 C_(H)=\frac{H(2H-1)\,\Gamma\left(H-\frac{1}{2}\right)}
{\Gamma\left(H+\frac{1}{2}\right)};
\end{equation}
see \eqref{eq:fbm-norm} on page \pageref{eq:fbm-norm}.
Also, $\delta=2H$, $\gamma=3-2H$ (see Table \ref{tab:rate}).
 Inequality (\ref{err-gen-overall}) becomes
\bel{err-fbm33}
\begin{split}
\sup_{0<t<T}\bE\|u(t)&-u_N^{n,r}(t)\|_{\bH}^2
\leq
C_H(T)\Bigg(
\frac{(C_1(H))^{N+1}T^{2H(N+1)}}{(N+1)!}
\bE\|u_0\|^2_{\bH}\\
&+
\frac{T^{2H}}{n^{2-2H}}\bE\|u_0\|^2_{\bH}+
\frac{T^{2H+2}}{n^{2-2H}}\bE\|u_0\|^2_{\bH^2}
\Bigg),
\end{split}
\ee
where $\limsup\limits_{T\to 0} C_H(T)>0$.
\end{example}

\subsection{Step-by-Step Approximation}

\subsubsection{Motivation}
\label{sec:ssa-motiv}

In the previous section, we constructed an
approximate solution for
equation (\ref{tr-evol1}) on the time interval
$[0,T]$ and
derived an error bound. The error bound suggests
that the
quality of the approximation improves for small
 values of $T$.
To construct the approximation for large values
 of $T$, it is natural to use
 a step-by-step method.

The main idea of the step-by-step method is as follows.
Let $\Psi_{t}$ be the solution operator for a
 homogenous linear  evolution equation,
 that is,
given an initial condition $u_0$, $u(t)=\Psi_tu_0$
is the solution
of the equation at time $t$. If the equation
is time-homogeneous
(has no explicit dependence on time,
such as (\ref{tr-evol1})) and
the solution is unique, then the solution operator has
the semi-group property:
\bel{sg-999}
\Psi_tu_0=\Psi_{t-s}u(s),\ t>s>0.
\ee
 If, for each $t> 0$, $u(t)$ is an element
of an infinite-dimensional Hilbert space $X$ with norm
$\|\cdot\|$, then a one-step approximate solution
$\ol{u}(t)$ can be constructed by
$$
\ol{u}(t)=\Pi^N\Psi_t u_0,
$$
where $\Pi^N$ is an orthogonal projection on
an $N$-dimensional  sub-space
of $X$. Assume that the approximation is of
order ${p/2}$ in time for some $p>1$:
\bel{aaa}
\|\ol{u}(t)-u(t)\|^2
= \|(I-\Pi^N)\Psi_tu_0\|^2 \leq Ct^p\|u_0\|^2,
\ee
where $I$ denotes the identity operator.
To construct a multi-step approximation on $[0,T]$, let
$0=t_0<t_1<\ldots<t_K=T$ be a uniform partition
of $[0,T]$ with step
$\tau$. Then define $\ol{u}_i$, $i=0,\ldots, K,$
 recursively as follows:
\begin{equation}
\ol{u}_0=u_0,\ \ol{u}_{i+1}=\Pi^N\Psi_{\tau} \ol{u}_i.
\end{equation}
For simplicity,
we assume that the initial condition is not
approximated, and concentrate only on
the effects of approximating the
solution operator $\Psi_t$. Then, by linearity,
\bel{bbb}
u(t_i)-\ol{u}_i=\Psi_{\tau} u(t_{i-1})
-\Pi^N\Psi_{\tau} \ol{u}_{i-1}
=(I-\Pi^N)\Psi_{\tau}u_{\tau_{i-1}}
+\Pi^N\Psi_{\tau}(u_{\tau_{i-1}}-\ol{u}_{i-1}).
\ee
By orthogonality, we find
\bel{ccc}
\begin{split}
\|u(t_i)-\ol{u}_i\|^2
&=\|\Psi_{\tau} u(t_{i-1})
-\Pi^N\Psi_{\tau} \ol{u}_{i-1}\|^2\\
&=\|(I-\Pi^N)\Psi_{\tau}u_{t_{i-1}}\|^2
+\|\Pi^N\Psi_{\tau}(u_{t_{i-1}}-\ol{u}_{i-1})\|^2.
\end{split}
\ee
Let $\Delta_i=\|u(t_i)-\ol{u}_i\|^2$.
Then (\ref{ccc}) and (\ref{aaa}) imply
\bel{ddd}
\Delta_i \leq C\tau^p \|u_{t_{i-1}}\|^2+
\|\Psi_{\tau} (u_{\tau_{i-1}}-\ol{u}_{i-1})\|^2.
\ee
In many situations, the semi-group $\Psi_t$ satisfies
\bel{sss}
\|\Psi_t f\|\leq e^{at} \|f\|
\ee
for some $a>0$. In this case,
 $\|u_{t_{i-1}}\|^2\leq e^{2aT}\|u_0\|^2$ and
(\ref{ddd}) implies
\bel{eee}
\Delta_i \leq C_1\tau^p\|u_0\|^2+e^{2a\tau}\Delta_{i-1},
\ee
or, after applying this inequality repeatedly,
\bel{eee1}
\Delta_i\leq C_1\tau^p\|u_0\|^2\sum_{j=0}^ie^{2a\tau j}.
\ee
Since
\begin{equation}
\sum_{j=0}^ie^{2a\tau j}
\leq \sum_{j=0}^Ke^{2a\tau j}
=\frac{e^{2a(K+1)\tau}}{e^{2a\tau}-1}
\leq \frac{e^{4aT}}{{e^{2a\tau}-1}},
\end{equation}
and $e^{2a\tau}-1\geq 2a\tau$, we conclude that
\begin{equation}
\Delta_i\leq C_2\tau^{p-1}\|u_0\|^2,
\end{equation}
that is,
\bel{eee2}
\max_{0\leq i\leq K}\|u(t_i)-\ol{u}_i\|^2
 \leq C_2\tau^{p-1}\|u_0\|^2,
\ee
where $C_2$ depends only on $T$ and
the semi-group $\Psi_t$.
In other words, the step-by-step approximation
 has order $(p-1)/2$ in time. The
derivation of this result essentially relies
on the following:
\begin{enumerate}
\item an approximation based on an orthogonal
 projection;
\item the property (\ref{sss}) of the
solution operator.
\end{enumerate}

\subsubsection{The Chaos Solution}
\label{sec:ssa-mr}

Let us consider equation (\ref{tr-evol1})
on  page \pageref{tr-evol1}.
The approximation $u_N^{n,r}(t)$
of the solution
is based on an orthogonal projection in the
space of square-integrable processes and,
 by Theorem \ref{th:chaos-m1},
page \pageref{th:chaos-m1}, the
solution operator for the equation
satisfies (\ref{sss}).  We
can therefore use (\ref{eee2}) to derive an error bound
 for the step-by-step
approximation of the solution of (\ref{tr-evol1}).

Let $0=t_0<t_1<\ldots<t_K=T$ be a uniform partition of
the interval $[0,T]$ with step $\tau$: $t_j=j\tau$,
 $j=0,\ldots, K$.
Let $\{m_k,\ k\geq 1\}$ be an orthonormal basis
in $L_2((0,T))$ and
$m_k^j(t)=m_k(t-t_j)(\chi_{t_{j+1}}(t)-\chi_{t_j}(t))$.
 We define
\begin{equation}
\xi_{k\ell}^j=\int_{t_{j-1}}^{t_j}m^j_k(t)dW_{\ell}(t),
\end{equation}
and then, for $\alpha\in \cJ$,
\begin{equation}
\xi_{\alpha}^j=\prod_{k,\ell}
\frac{H_{\alpha_{k\ell}}(\xi_{k\ell}^j)}
{\sqrt{\alpha_{k\ell}!}}.
\end{equation}
Note that the random variables
$\xi_{k\ell}^i$ and $\xi_{pq}^j$ are independent
for different $i,j$.

If $u=u(t;u_0)$ is the square-integrable chaos
solution of the homogeneous equation (\ref{tr-evol1})
with initial condition $u_0$,
then, by uniqueness and time homogeneity, we have
\bel{sg-999-1}
u(t_j;u_0)=u(\tau;u(t_{j-1},u_0)),
\ee
which is a particular case of the general
relation (\ref{sg-999}).
Also, by Theorem \ref{th:chaos-m1}
on page \pageref{th:chaos-m1},
\bel{sss-999}
\bE\|u(t;u_0)\|^2_{\bH}
\leq e^{(C_A+C_B)t}\bE \|u_0\|_{\bH}^2,
\ee
 which is a particular case
of (\ref{sss}) on page \pageref{sss}.

Next, consider the following modification
of the S-system (\ref{err-s-syst})
 from page \pageref{err-s-syst}:
\bel{err-s-syst99}
 \begin{split}
&u^j_{(0)}(t;f)=f+\int_{t_{j-1}}^t \cA u_{(0)}(s;f)ds,
\ \ |\alpha|=0;\\
& u^j_{\alpha}(t;f)=
\int_{t_{j-1}}^t\cA u^j_{\alpha}(s;f) ds
+
\sum_{k,\ell=1}^{\infty}\sqrt{\alpha_{k\ell}}
\int_{t_{j-1}}^t\cB_{\ell}
u^j_{\alpha-\epsilon_{k\ell}}(s;f)\, (\cK m_k^j)(s)ds
\end{split}
 \ee
 for $\ |\alpha|>0$, where $\ t_{j-1}\leq t\leq t_j,$
  $\bE\|f\|_{\bH}^2<\infty$, and
$f$ is random but independent of  $\xi_{\alpha}^j,\
\alpha \in \cJ$.

\begin{theorem}
\label{th:rec-S}
Assume {\bf A1}--{\bf A5}. Then, for $j=1,\ldots, K$,
\bel{eq:wc-rec}
 u(t_j;u_0)
 =\sum_{\alpha\in \cJ}
 u^{j}_{\alpha}(t_j;u(t_{j-1};u_0))\xi_{\alpha}^j
 \ee
 and
 \bel{eq:wc-rec1}
 \bE \|u(t_j;u_0)\|_{\bH}^2
 \leq e^{(C_A+C_B)\tau}\bE\|u(t_{j-1};u_0)\|_{\bH}^2.
 \ee
 \end{theorem}

 \begin{proof}
 Since  $u(t_{j-1};u_0)$ is independent
 of $\{\xi_{\alpha}^j,\ \alpha\in \cJ\},$
  Theorem \ref{th:chaos-m1} can be applied
  on each interval $[t_{j-1},t_j]$,  $j=1,\ldots,K,$
  with $u_0$ replaced by $u(t_{j-1};u_0)$.
  Then (\ref{sg-999-1}) becomes \eqref{eq:wc-rec}
  and \eqref{eq:reg-sol} becomes \eqref{eq:wc-rec1}.
 \end{proof}

Define the multi-step approximation of $u$ as follows:
\bel{multistep00}
u^{n,r}_N(t_j)=\sum_{\alpha\in \cJ^{n,r}_N}
u^{j}_{\alpha}(t_j;u^{n,r}_N(t_{j-1}))\xi_{\alpha}^j,
\ j=1,\ldots, K,
\ee
with $u^{n,r}_N(t_0)=u_0$.
The following theorem provides an error bound for this
approximation.

\begin{theorem}\label{th:multipstep}
Under the assumptions of Theorems \ref{th:err-b-2}
and \ref{th:err-b-3}
\bel{step-by-step-err}
\begin{split}
\max_{j=1,\ldots,K}
\bE\|u(t_j;u_0)&-u_N^{n,r}(t_j)\|_{\bH}^2 \leq
C(T)\Bigg(\frac{(\tau C_B)^{N}}{(N+1)!}
\bE\|u_0\|^2_{\bH}\\
&+\frac{\tau^{\delta_1-1}}{n^{\gamma_1-1}}
\bE\|u_0\|^2_{\bH}
+\frac{\tau^{\delta}+1}{n^{\gamma-1}}
\bE\|u_0\|^2_{\bH^2}+
\varepsilon(r)\bE\|u_0\|^2_{\bH}\Bigg).
\end{split}
\ee
\end{theorem}

\begin{proof} Since $u^{n,r}_N(\tau)$ is an orthogonal
projection of $u(\tau;u_0)$
on the span of $\xi_{\alpha}, \ \alpha \in \cJ^{n,r}_N$,
 and (\ref{sss-999}) holds, the result follows
 from (\ref{err-gen-overall1}) and (\ref{eee2}).
\end{proof}

\begin{example}
Consider the equation
$$
u(t)=u_0+\int_0^tu_{xx}(s)ds+\mfX(\chi_tu),\
0\leq t\leq T,\ x\in \bR.
$$
With only one noise driving the equation, we
have $u_N^{n,r}=u_N^{n}$. Also,
$\bH^{r}=\bH^{r}(\bR)$
is the Sobolev space \eqref{sobolev} and
$\bH=L_2(\bR)$.

(a) If $\mfX$ is an
 {Ornstein-Uhlenbeck noise}  with
 parameter $b$, then
$C_B=(1+\sqrt{b\tau})^2,$
 $\delta=1$, $\gamma=2$, $\delta_1=3$, $\gamma_1=4$
 (see Table \ref{tab:rate}),
  so that (\ref{step-by-step-err}) becomes
\bel{err-OU-ss}
\begin{split}
\max_{j=1,\ldots,K}
\bE\|u(t_j;u_0)-u_N^{n,r}(t_j)\|_{\bH}^2 &\leq
C_b(T)\Bigg(
\frac{(1+\sqrt{b\tau})^{2N+2}\,\tau^{N}}{(N+1)!}
\bE\|u_0\|^2_{\bH}\\
&+
\frac{\tau^{2}}{n^{3}}\bE\|u_0\|^2_{\bH}+
\frac{\tau^2}{n}\bE\|u_0\|^2_{\bH^2}\Bigg).
\end{split}
\ee

(b) If $\mfX$ is an
 $H$-fractional white noises with
parameter $H\in (1/2,1)$, then
 $C_B=C_1(H){\tau}^{2H-1}$ (see \eqref{eq:C1H}),
$\delta=\delta_1=2H$, $\gamma=\gamma_1=3-2H$
(see Table \ref{tab:rate}),
so that (\ref{step-by-step-err}) becomes
\bel{err-fbm33-ss}
\begin{split}
\max_{j=1,\ldots,K}
\bE\|u(t_j;u_0)-u_N^{n,r}(t_j)\|_{\bH}^2 & \leq
C_H(T)\Bigg(
\frac{(C_1(H))^{N+1}\,\tau^{2H(N+1)-1}}{(N+1)!}
\bE\|u_0\|^2_{\bH}\\
&+\frac{\tau^{2H-1}}{n^{2-2H}}\bE\|u_0\|^2_{\bH}+
\frac{\tau^{2H+1}}{n^{2-2H}}\bE\|u_0\|^2_{\bH^2}
\Bigg).
\end{split}
\ee
\end{example}

\section*{Acknowledgment}
SVL is grateful to Professor Wendell Fleming
for very fruitful discussions, and to the
Division of Applied Mathematics at Brown University  for
 hospitality and stimulating atmosphere.


\begin{thebibliography}{10}

\bibitem{AMN1}
E.~Al\`{o}s, O.~Mazet, and D.~Nualart, \emph{{Stochastic Calculus With Respect
  to Fractional Brownian Motion with Hurst Parameter Less Than $\frac{1}{2}$}},
  Stochastic Process. Appl. \textbf{86} (2000), no.~1, 121--139.

\bibitem{AMN}
\bysame, \emph{{Stochastic Calculus With Respect to Gaussian Processes}}, Ann.
  Probab. \textbf{29} (2001), no.~2, 766--801.

\bibitem{A}
A.~Amirdjanova, \emph{{Nonlinear Filtering with Fractional Brownian Motion}},
  Appl. Math. Optim. \textbf{46} (2002), no.~2--3, 81--88.

\bibitem{AC}
A.~Amirdjanova and S.~Chivoret, \emph{New method for optimal nonlinear
  filtering of noisy observations by multiple stochastic fractional integral
  expansions}, Comput. Math. Appl. \textbf{52} (2006), no.~1-2, 161--178.
  \MR{MR2262164}

\bibitem{BP}
R.~J. Barton and H.~V. Poor, \emph{{Signal Detenction in Fractional Gaussian
  Noise}}, IEEE Trans. Infom. Theory \textbf{34} (1988), no.~5,part 1,
  943--959.

\bibitem{Bennaton}
J.~F. Bennaton, \emph{{Discrete Time {Galerkin} Approximation to the Nonlinear
  Filtering Solution}}, J. Math. Anal. Appl. \textbf{110} (1985), no.~2,
  364--383.

\bibitem{BGR}
A.~Bensoussan, R.~Glowinski, and R.~Rascanu, \emph{{Approximations of the
  {Zakai} Equation by the Splitting-up Method}}, SIAM J. Control Optim.
  \textbf{28} (1990), no.~6, 1420--1431.

\bibitem{BH}
T.~Bj\"{o}rk and H.~Hult, \emph{{A Note on Wick Products and the Fractional
  Black-Scholes Model}}, Finance Stoch. \textbf{9} (2005), no.~2, 197--209.

\bibitem{LeB}
A.~Le Breton, \emph{{Filtering and Parameter Estimation in a Simple Linear
  System Driven by Fractional Brownian Motion}}, Statist. Probab. Lett.
  \textbf{38} (1998), no.~3, 263--274.

\bibitem{BdKl}
A.~Budhiraja and G.~Kallianpur, \emph{{Approximations to the Solution of the
  {Zakai} Equation Using Multiple {Wiener} and {Stratonovich} Integral
  Expansions}}, Stochastics Stochastics Rep. \textbf{56} (1996), no.~3--4,
  271--315.

\bibitem{CM}
R.~H. Cameron and W.~T. Martin, \emph{{The Orthogonal Development of Nonlinear
  Functionals in Series of {Fourier-Hermite} Functionals}}, Ann. of Math.
  \textbf{48} (1947), no.~2, 385--392.

\bibitem{DH}
W.~Dai and C.~C. Heyde, \emph{{It\^{o}'s Formula with Respect to Fractional
  Brownian Motion and its Application}}, J. Appl. Math. Stochastic Anal.
  \textbf{9} (1996), no.~4, 439--448.

\bibitem{DU}
L.~Decreusefond and A.~S. \"{U}st\"{u}nel, \emph{{Stochasic Analysis of the
  Fractional Brownian Motion}}, Potential Anal. \textbf{10} (1998), no.~2,
  177--214.

\bibitem{DHP}
T.~E. Duncan, Y.~Hu, and B.~Pasik-Duncan, \emph{{Stochastic Calculus for
  Fractional Brownian Motion I: Theory}}, SIAM J. Control Optim. \textbf{38}
  (2000), no.~2, 582--612.

\bibitem{GV}
I.~M. Gelfand and N.~J. Vilenkin, \emph{{Generalized Functions IV: Applications
  of Harmonic Analysis}}, Acadimic Press, 1964.

\bibitem{GermPicc}
A.~Germani and M.~Piccioni, \emph{Semi-discretisation of stochastic partial
  differential equations on {${\bf R}^d$} by a finite element technique},
  Stochastics \textbf{23} (1988), no.~2, 131--148.

\bibitem{HKPS}
T.~Hida, H-H. Kuo, J.~Potthoff, and L.~Sreit, \emph{{White Noise}}, Kluwer,
  1993.

\bibitem{HOUZ}
H.~Holden, B.~{\O}ksendal, J.~Ub{\o}e, and T.~Zhang, \emph{{Stochastic Partial
  Differential Equations: A Modeling, White Noise Functional Approach}},
  Birkh\"{a}user, 1996.

\bibitem{HO}
Y.~Hu and B.~{\O}ksendal, \emph{{Fractional White Noise Calculus and
  Applications to Finance}}, Infin. Dimens. Anal. Quantum Probab. Relat. Top.
  \textbf{6} (2003), no.~1, 1--32.

\bibitem{KIto}
K.~Ito, \emph{{Approximation of the {Zakai} Equation for Nonlinear Filtering}},
  SIAM J.~ Control~Optim. \textbf{34} (1996), no.~2, 620--634.

\bibitem{ItoRoz}
K.~Ito and B.~L. Rozovskii, \emph{{Approximation of the {Kushner} equation for
  nonlinear filtering}}, SIAM J.~Control~Optim. \textbf{38} (2000), no.~3,
  893--915.

\bibitem{KBR}
M.~L. Kleptsyna, A.~{Le} Breton, and M.-C. Roubaud, \emph{{General Approach to
  Filtering with Fractional Brownian Noises: Application to Linear Systems}},
  Stochastics Stochastics Rep. \textbf{71} (2000), no.~1--2, 119--140.

\bibitem{KKA}
M.~L. Kleptsyna, P.~E. Kloeden, and V.~V. Anh, \emph{{Existence and Uniqueness
  Theorems for Stochastic Differential Equations with Fractal Brownian
  Motion}}, Problems Inform. Transmission \textbf{34} (1998), no.~4, 332--341.

\bibitem{KlPl}
P.~E. Kloeden and E.~Platen, \emph{{Numerical Solution of Stochastic
  Differential Equations}}, Springer, 1992.

\bibitem{KPS}
S.~G. Krein, Ju.~I. Petunin, and E.~M. Semeonov, \emph{{Interpolation of Linear
  Operators}}, AMS, 1982.

\bibitem{Lin}
S.~J. Lin, \emph{{Stochastic Analysis of Fractional Brownian Motions}},
  Stochastics Stochastics Rep. \textbf{55} (1995), no.~1--2, 121--140.

\bibitem{LMR}
S.~V. Lototsky, R.~Mikulevicius, and B.~L. Rozovskii, \emph{{Nonlinear
  Filtering Revisited: A Spectral Approach}}, SIAM J. Contr. Optim. \textbf{35}
  (1997), no.~2, 435--461.

\bibitem{LR1}
S.~V. Lototsky and B.~L. Rozovskii, \emph{{Wiener Chaos Solutions of Linear
  Stochastic Evolution Equations}}, Ann. Probab. \textbf{34} (2006), no.~2,
  638--662.

\bibitem{LtSt}
S.~V. Lototsky and K.~Stemmann, \emph{{From Random Processes to Generalized
  Fields: a Unified Approach to Stochastic Integration}}, Submitted to {\em
  Stochastic Processes and Their Applications}.

\bibitem{Maj}
P.~Major, \emph{{Multiple Wiener-It{\^{o}} Integrals. With Applications to
  Limit Theorems}}, {Lecture Notes in Mathematics}, vol. 849, Springer, 1981.

\bibitem{Mal}
P.~Malliavin, \emph{Stochastic analysis}, Springer, 1997.

\bibitem{Mil}
G.~Milstein, \emph{{Numerical Integration of Stochastic Differential
  Equations}}, Kluwer, 1995.

\bibitem{Nualart2}
D.~Nualart, \emph{{Stochastic Integration With Respect to Fractional Brownian
  Motion and Applications}}, Stochastic models (Mexico City, 2002), Contemp.
  Math., vol. 336, Amer. Math. Soc., Providence, RI, 2003, pp.~3--39.

\bibitem{Nualart}
D.~Nualart, \emph{The {M}alliavin calculus and related topics}, second ed.,
  Probability and its Applications (New York), Springer-Verlag, Berlin, 2006.
  \MR{MR2200233 (2006j:60004)}

\bibitem{Oksendal}
B.~{\O}ksendal, \emph{Stochastic differential equations}, sixth ed.,
  Universitext, Springer-Verlag, Berlin, 2003. \MR{MR2001996 (2004e:60102)}

\bibitem{PT}
V.~Pipiras and M.~S. Taqqu, \emph{{Integration Questions Related to Fractional
  Brownian Motion}}, Probab. Theory Related Fields \textbf{118} (2000), no.~2,
  251--291.

\bibitem{Roz}
B.~L. Rozovskii, \emph{{Stochastic Evolution Systems}}, Kluwer, 1990.

\bibitem{Za}
M.~Z\"{a}hle, \emph{{Integration With Respect to Fractal Functions and
  Stochastic Calculus I}}, Probab. Theory Related Fields \textbf{111} (1998),
  no.~3, 333--374.

\end{thebibliography}

\providecommand{\bysame}{\leavevmode\hbox to3em{\hrulefill}\thinspace}
\providecommand{\MR}{\relax\ifhmode\unskip\space\fi MR }
\providecommand{\MRhref}[2]{%
  \href{http://www.ams.org/mathscinet-getitem?mr=#1}{#2}
}
\providecommand{\href}[2]{#2}

\end{document}